\documentclass[12pt, amsfonts]{article}

\usepackage[latin1]{inputenc}
\usepackage[all]{xy}

\usepackage{amsfonts,amssymb,amsmath,latexsym}

\begin{document}
\pretolerance 1000

\centerline
{\it \Large     On the Structure of the  Solution Set  of a Sign }
\vskip.1cm

\centerline
{\it \Large    Changing  Perturbation of the p-Laplacian  }
\vskip.1cm

\centerline
{\it \Large   under Dirichlet  Boundary Condition\footnote{Supported   by  CNPq/CAPES/PROCAD/UFG/UnB-Brazil}}

\numberwithin{equation}{section}
\newtheorem{theorem}{Theorem}[section]
\newtheorem{lemma}{Lemma}[section]
\newtheorem{remark}{Remark}[section]
\newtheorem{corollary}[theorem]{Corollary}
\newtheorem{proposition}[theorem]{Proposi\c c\~ao}
\newtheorem{definition}{Definition}[section]

\newcommand{\prova}{\normalfont{\bf Proof }}
\newcommand{\R}{I\!\!R}
\newcommand{\rn}{I\!\!R^{N}}
\newcommand{\N}{I\!\!N}
\newcommand{\Z}{I\!\!Z}

\newcommand{\C}{\mathbb{C}}
\newcommand{\w}{W_{0}^{1,p}(\Omega)}
\newcommand{\nd}{\noindent}
\bigskip

\begin{center}
{J. V.  Goncalves~~ M. R. Marcial }\\
\smallskip
  \scriptsize{Universidade Federal de Goi\'as\\
   Instituto de Matem\'atica e Estat\'istica\\
   74001-970 Goi\^ania, GO - Brasil}\\
Email: goncalves.jva@gmail.com   
\end{center}

\medskip
\begin{abstract}
\vskip.3cm

\nd In a recent  paper D. D. Hai showed that  the equation  $ -\Delta_{p} u = \lambda f(u)~\mbox{in}~ \Omega$, under Dirichlet boundary condition, where $\Omega \subset {\bf R^N}$ is a bounded  domain with smooth boundary $\partial\Omega$, $\Delta_{p}$ is the   p-Laplacian, 
 $f : (0,\infty) \rightarrow {\bf R} $ is a continuous function which may blow up to  $\pm \infty$ at the origin, admits a solution if $\lambda > \lambda_0$ and has no solution if $0 < \lambda < \lambda_0$. In this paper we show that  the solution set $\mathcal{ S}$ of the equation above, which is not empty by Hai's results, actually  admits  a continuum  of positive solutions. 
\vskip.5cm

\end{abstract}
\vskip.2cm

\section{Introduction}

\nd In this paper we establish existence of a continuum of positive solutions of 
$$
\left\{\begin{array}{rllr}
-\Delta_{p} u&=&\lambda f(u) +h~~  \mbox{in}~~  \Omega,\\

u&=&0~~ \mbox{on}~~  \partial {\Omega,}
 \end{array}
\right.\leqno{(P)_\lambda}
$$
\nd where  $\Omega \subset {\bf R^N}$ is a bounded domain with smooth boundary   $\partial {\Omega}$,   $\lambda>0$ is a real parameter, $f : (0,\infty) \rightarrow {\bf R} $ is a continuous function  which may blow up to  $\pm \infty$ at the origin and  $h : \Omega \rightarrow {\bf R}$  is a nonnegative $L^\infty$-function.
\vskip.2cm

\begin{definition}
\nd By a solution of $(P)_\lambda$ we mean a function $u\in \w$ such that
  \begin{equation}\label{eq 0.1}
  \int_\Omega|\nabla u|^{p-2}\nabla u.\nabla\varphi dx=\lambda\int_\Omega f(u)\varphi dx+\int_\Omega h\varphi dx,~~ \varphi\in\w.
  \end{equation}
 \end{definition}
\vskip.2cm

\begin{definition}
\nd The solution set of ${(P)_\lambda}$ is 
  \begin{equation}\label{eq 0.2}
  \mathcal{ S}:=\big{\{} (\lambda,u)\in(0,\infty)\times C(\overline{\Omega})\ \big{|}\  \mbox{u is a solution of}~ (P)_\lambda \  \big{\}}.
  \end{equation}
 \end{definition}
\vskip.2cm

\nd It was shown by Hai \cite{H} that there is a positive number $\lambda_0$ such that ${(P)_\lambda}$ admits:  a solution if $\lambda > \lambda_0$ and  no solution if $\lambda < \lambda_0$. Our aim is to investigate existence of connected components of  $\mathcal{ S}$. By adapting estimates in \cite{H} we succeeded in showing the existence of a continuum $ \Sigma \subset \mathcal{ S}$ such that  $Proj_{\bf R } \Sigma = (\lambda_0, \infty)$.
\vskip.2cm

\nd The assumptions on $f$ are:
\begin{itemize}
 \item[$(f)_{1}$]~~  $f:(0,\infty) \rightarrow {\bf R}$ is continuous and
 $$
\displaystyle \lim_{u \to \infty}\frac{f(u)}{u^{p-1}}=0,
$$
 \item[$(f)_{2}$]~~  there are positive numbers $a, \beta, A$ with $\beta<1$ such that
$$
\mbox{\rm (i)}~~   f(u) \geq\frac{a}{u^{\beta}}~ \mbox{for}~u>A,~~~
\mbox{\rm (ii)}~~  \displaystyle \limsup_{u \to 0} u^{\beta}|f(u)|< \infty.
 $$
 \end{itemize}

\nd We give below a few examples of functions $f$ satisfying  $(f)_1$, $ (f)_2$. Those functions appear in several earlier works on existence of solutions, cf. section 2.,
\vskip.3cm

  {\bf a)}~~ $\displaystyle u^{q}-\frac{1}{u^{\beta}}$,~ $\beta>0$,~ $0 < q<p-1$,~~~ 
  {\bf b)}~ ~ $\displaystyle \frac{1}{u^{\beta}}-\frac{1}{u^{\alpha}}$,~~ $0 < \beta < \alpha < 1$,\\
 
{\bf c)}~ ~ $\displaystyle a-\frac{1}{u^{\alpha}}$, ~~ $a>0,~ 0< \alpha<1$,~~~ {\bf d)} $\displaystyle \frac{1}{u^{\alpha}} + u^q$,~ $ 0 < \alpha<1,~0 < q < p-1$,\\

 {\bf e)}~~ $\displaystyle \frac{1}{u^{\alpha}}$, $0 < \alpha<1$,~~ {\bf f)}~ $\displaystyle  \ln u$.

\vskip.5cm
\nd The main results of this paper are,

\begin{theorem}\label{teo 0.1}
  Assume  $(f)_1-(f)_2$. Then there is a number $\lambda_*>0$ and a connected subset
$ \Sigma$ of~  $ [\lambda_*,\infty) \times C(\overline{\Omega})$
 satisfying,
  \begin{equation}\label{eq 1.1}
  {\Sigma  } \subset  \mathcal{ S},
  \end{equation}
  \begin{equation}\label{eq 1.2}
 \ {\Sigma } \cap\big(\{\lambda \}\times C(\overline{\Omega})\big)\neq \emptyset,~~ \lambda_* \leq \lambda < \infty.
  \end{equation}
\end{theorem}
\vskip.2cm

\nd The prove of theorem  \ref{teo 0.1} will be achieved by at first proving the following result.
\vskip.2cm

\begin{theorem}\label{teo 0.0}
  Assume $(f)_1-(f)_2$. Then there is a number $\lambda_*>0$ and for each  $\Lambda>\lambda_*$ there is a connceted set   $\Sigma_{\Lambda}\subset ([\lambda_*,\Lambda]\times C(\overline{\Omega})$ satisfying
  \begin{equation}\label{eq 1.1}
  \Sigma_{\Lambda}\subset \mathcal{\huge S},
  \end{equation}
  \begin{equation}\label{eq 1.2}
  \Sigma_{\Lambda}\cap\big(\{\lambda_*\}\times C(\overline{\Omega})\big)\neq \emptyset,
  \end{equation}
  \begin{equation}\label{eq 1.3}
  \Sigma_{\Lambda}\cap\big(\{\Lambda\}\times C(\overline{\Omega})\big)\neq \emptyset.
  \end{equation}
\end{theorem}
\vskip.2cm

\begin{remark}
 The present work  is motivated by  Hai \cite{H}. We will use  $C, C_1, C_2,  \widetilde{C}$ to denote positive cumulative constants.
\end{remark}

\section{Background}

\nd The Dirichlet problem
\begin{equation}\label{gen_p_probl}
-\Delta_{p} u =  f(x,u)~~ \mbox{in}~~ \Omega,~~
u = 0~ \mbox{on}~ \partial {\Omega,} 
\end{equation}
\nd where $f: \Omega \times (0,\infty) \to {\bf R}$ is a function  satisfying a condition like  $f(x,r) \to +\infty$ as $r \to 0$, referred to as singular at the origin has been extensively studied in the last years.
\vskip.2cm

\nd In the pioneering work  \cite{CRT}, it was shown by Crandall, Rabinowitz \& Tartar through the use of topological methods, e.g. Schauder Theory and Maximum Principles, that the problem 
$$
\left\{\begin{array}{rllr}
-\Delta u&=& u^{-\gamma}~~ \mbox{in}~~ \Omega,\\
u &> & 0~~ \mbox{in}~~ \Omega,\\
u&=&0~~ \mbox{on}~~ \partial {\Omega},
 \end{array}
\right.
$$
\nd where   $\gamma>0$, admits a solution  $u \in C^2(\Omega)\cap C(\overline{\Omega})$, (see also the references of \cite{CRT}).
\vskip.2cm

\nd Subsequently,  Lazer \& McKenna in \cite{LM}, established, among other results,  the existence of a solution $u\in C^{2+\alpha}(\Omega)\cap C(\overline{\Omega})$  ($0 < \alpha < 1$)   for the problem
$$
\left\{\begin{array}{rllr}
-\Delta u&=& p(x)u^{-\gamma}~~ \mbox{in}~~  \Omega,\\
u &> & 0~~  \mbox{in}~~  \Omega,\\
u&=&0~~  \mbox{on}~~  \partial {\Omega},
 \end{array}
\right.
$$
\nd where   $p\in C^\alpha(\overline{\Omega})$ is a positive function.
\vskip.2cm

\nd Several techniques have been employed in the study of (\ref{gen_p_probl}).  In  \cite{Z}, by using lower and upper solutions,  Zhang showed that there is  some number  $\overline{\lambda}\in(0,+\infty)$    such that the problem
 $$
\left\{\begin{array}{rllr}
\displaystyle-\Delta u+\frac{1}{u^{\alpha}} &=& \lambda u^p~~  \mbox{in}~~  \Omega,\\
u &> & 0~~  \mbox{in}~~  \Omega,\\
u&=&0~~ \mbox{in}~~  \partial {\Omega,} 
\end{array}
\right.
$$
\nd where   $\alpha, p \in(0,1)$, admits a solution  $u_\lambda\in C^{2+\gamma}(\Omega)\cap C(\overline{\Omega})\cap H_0^1(\Omega)$
with $u^{-\alpha}_{\lambda} \in L^{1}(\Omega)$ for each  $\lambda>\overline{\lambda}$ and  no solution in $ C^{2}(\Omega)\cap C(\overline{\Omega})$ for  $\lambda<\overline{\lambda}$. It was also shown that the problem above admits no solution in $C(\overline{\Omega})\cap H_0^1(\Omega)$ if
 $\alpha\geq 1$, $\lambda>0$  and $p>0$.
\vskip.2cm

\nd In  \cite{JST},  Giacomoni,  Schindler \& Takac employed variational methods to  investigate the problem
$$
\left\{\begin{array}{rllr}
\displaystyle -\Delta_{p} u &=&  \frac{\lambda}{u^\delta}+u^q ~ \mbox{in}~  \Omega,\\
 u&>0&~~  \mbox{in}~~  \Omega,\\
u &= 0&~~ \mbox{on}~~    \partial\Omega,
\end{array}
\right.
$$
\nd where  $1<p<\infty$, $p-1<q<p*-1$, $\lambda>0$ and $0<\delta<1$ with $p^*=\frac{Np}{n-p}$ if $1<p<N$, $p^*\in(0,\infty)$ large  if $p=N$ and $p^*=\infty$ if $p>N$.  Several results were shown in that paper, among them  existence, multiplicity and regularity of solutions.
\vskip.2cm

\nd In  \cite{PZ}, Perera \& Zhang  used variational methods to prove existence of solution for the problem
$$
\left\{\begin{array}{rllr}
-\Delta_{p} u&=& a(x)u^{-\gamma}+\lambda f(x,u)~~  \mbox{in}~~ \Omega,\\
u&>& 0~~  \mbox{in}~~  {\Omega},\\
u&=&0~~  \mbox{on}~~  \partial {\Omega},
\end{array}
\right.
$$
\nd where   $1<p<\infty$, $\gamma, \lambda>0$ are numbers, $a\geq 0$ is  a measurable, not identically zero function and $f:\Omega\times[0,\infty)\to\R$ is a  Carath\'eodory  satisfying
$$
\sup_{(x,t)\in\Omega\times[0,T]} |f(x,t)|<\infty
$$
\nd for each $T > 0$.
\vskip.2cm

\nd There is a broad literature on singular problems and we further refer the reader to 
Gerghu \& Radulescu  \cite{GR},  Goncalves, Rezende \&  Santos \cite{JVA},  Hai \cite{H-1, H}, Mohammed  \cite{M}, Shi \& Yao  \cite{SY},  Hoang Loc \& Schmitt \cite{LS}, Montenegro \& Queiroz \cite{MQ}   and their references. 

\section{ Some Auxiliary  Results }

\nd We gather below a few technical results. For completeness, a few  proofs will be provided in the Appendix.  The Euclidean distance from $x \in \Omega$ to $\partial \Omega$ is
$$
d(x) = dist(x,\partial\Omega).
$$
\nd The result below derives from   Gilbarg $\&$ Trudinger \cite{GT}, V\`azquez \cite{V}.
\vskip.5cm

\begin{lemma}\label{lem 3.2}
\nd Let $\Omega \subset {\bf R}^{N}$ be a  smooth,  bounded, domain.  Then
\begin{itemize}
\item[$\mbox{\rm (i)}$] $d \in Lip(\overline{\Omega})$ and  $d$ is $C^2$ in a neighborhood of  $\partial\Omega$,
\item[$\mbox{\rm (ii)}$] if  $\phi_1$ denotes a positive eigenfunction of  $(-\Delta_p, \w)$ one has,
$$
\phi_1 \in C^{1,\alpha}(\overline{\Omega})~\mbox{with}~~0  < \alpha < 1,~ \frac{\partial \phi_1}{\partial \nu} < 0~\mbox{on}~ \partial \Omega,
$$
\nd and there are positive constants $C_1, C_2$  such that
$$
 C_1 d(x) \leq \phi_1(x) \leq C_2 d(x),~  x\in\Omega.
$$
\end{itemize}
\end{lemma}

\nd The result below is due to  Crandall, Rabinowitz $\&$ Tartar $\cite{CRT}$, Lazer $\&$ McKenna $\cite{LM}$ in the case $p =2$ and
Giacomoni, Schindler $\&$ Takac \cite{JST} in the  case
$1<p<\infty$.

\begin{lemma}\label{lem 3.1}
Let $\beta \in (0,1)$ and  $m>0$. Then the problem
\begin{equation}
\left\{\begin{array}{rcll}\label{eq 1.1}
-\Delta_p u &=& \frac{m}{u^\beta} & \mbox{in}\  \Omega,\\
u &>&0 &\mbox{in}\  \Omega,\\
 u  &=& 0                                 &\mbox{on}\  \partial\Omega,
\end{array}\right.
\end{equation}
\nd admits an only weak solution  $u_m \in \w$. Moreover
$u_m\geq \epsilon_m\phi_1~ \mbox{in}~ \Omega$ for some constant $\epsilon_m > 0$.
\end{lemma}

\begin{remark}\label{on lem 3.2}
 \nd  By the results in   $\cite{L, JST}$, there is $\alpha \in(0,1)$  such that $
u_m\in C^{1,\alpha}(\overline{\Omega})$.
\end{remark}
\vskip.1cm

\nd The  result below, which is crucial in this work, and whose proof is provided in the Appendix, is basically due to  Hai \cite{H}.
\begin{lemma}\label{lem 3.3}
Let  $g\in L_{loc}^{\infty}(\Omega)$. Assume  that there is $\beta \in(0,1)$  and  $C>0$ such that
\begin{equation} \label{eq 3.7}
|g(x)|\leq \frac{C}{d(x)^{\beta}},~x \in \Omega.
\end{equation}
\nd Then there is an only  weak solution   $u\in\w$  of
\begin{equation} \label{eq 3.8}
\left\{
\begin{array}{rclr}
-\Delta_{p} u&=& g& \mbox{in}~~ \Omega\\
u&=&0&  \mbox{on} ~ \partial {\Omega} .
\end{array}
\right.
\end{equation}
\nd In addition,  there exist constants  $\alpha\in(0,1)$ and $M>0$,  with $M$ depending only on   $C, \beta, \Omega$ such that  $u\in C^{1,\alpha}(\overline{\Omega})$ and $||u||_{C^{1,\alpha}(\overline{\Omega})} \leq M$.
\end{lemma}

\begin{remark}\label{SOL OPER}
\nd The solution operator associated to $(\ref{eq 3.8})$ is: let
$$
{{\mathcal{M}}_{ \beta, \infty} }= \{ g \in L^{\infty}_{loc}(\Omega)~|~ |g(x) | \leq \frac{C}{d(x)^{\beta} },~x \in \Omega \},
$$

$$
S : {{\mathcal{M}}_{ \beta, \infty} } \rightarrow  \w \cap C^{1,\alpha}(\overline{\Omega}),~~ S(g) := u.
$$
\nd  Notice that
$$
||S(g)||_{C^{1,\alpha}(\overline{\Omega})} \leq M,
$$
\nd for all  $g\in {{\mathcal{M}}_{C,d, \beta, \infty} }$  with $M$ depending only on   $C, \beta, \Omega$.
\end{remark}

\begin{corollary}\label{cor 3.2}
Let  $g, \widetilde{g}\in L^\infty_{loc}(\Omega)$ with $g\geq 0$, $g\ne 0$ satisfying $(\ref{eq 3.7})$. Then, for each $\epsilon>0$,  the problem
\begin{equation}\label{eq 3.11}
\left\{\begin{array}{rlllr}
-\Delta_p u_\epsilon & = & g~ \chi_{\{d >\epsilon\}}+\widetilde{g}~ \chi_{\{d <\epsilon\}} & \mbox{em} & \Omega;\\
          u_\epsilon & = & 0                                                       & \mbox{em} & \partial\Omega
\end{array}\right.
\end{equation}
\nd admits an only solution  $u_\epsilon \in C^{1,\alpha}(\overline{\Omega})$ for some $\alpha\in(0,1)$. In addition, there is  $\epsilon_0>0$ such that
$$
u_\epsilon\geq \frac{u}{2} \quad \mbox{in}\quad \Omega \quad \mbox{for each}\quad \epsilon\in(0,\epsilon_0),
$$
\nd where $u$ is the solution of  $(\ref{eq 3.8})$.
\end{corollary}
\vskip.1cm

\nd A proof of the Corollary above will be included in the  Appendix.

\section{Existence of  Lower and Upper Solutions }

\nd In this section we present two results, essentially due to  Hai \cite{H}, on existence of lower and upper solutions of $(P)_\lambda$. At first some definitions.
\vskip.2cm

\begin{definition}
A function  $\underline{u} \in \w$ with $\underline{u}>0~\mbox{in}~\Omega$ such that
$$
\int_\Omega|\nabla \underline{u}|^{p-2}\nabla \underline{u}.\nabla\varphi dx\leq\lambda\int_\Omega f(\underline{u})\varphi dx+\int_\Omega h\varphi dx,~ \varphi\in\w,~\varphi\geq 0
$$
 is a lower solution of $(P)_\lambda$.
\end{definition}
  \begin{definition}
A function   $\overline{u} \in \w$ with $\overline{u}>0~\mbox{in}~\Omega$ such that
$$
\int_\Omega|\nabla \overline{u}|^{p-2}\nabla \overline{u}.\nabla\varphi dx\geq\lambda\int_\Omega f(\overline{u})\varphi dx+\int_\Omega h\varphi dx, ~ \varphi\in\w,\ \varphi\geq 0.
$$
 is an upper  solution of $(P)_\lambda$.
 \end{definition}

\begin{theorem} \label{teo 0.2}
Assume $(f)_1-(f)_2$. Then there  exist  $\lambda_*>0$ and a non-negative  function $\psi \in C^{1, \alpha}(\overline{\Omega})$,  with $\psi > 0~  \mbox{in}~  {\Omega}$, 
$ \psi = 0~  \mbox{on}~ \partial {\Omega}$,  $\alpha \in (0,1)$ such that for each  $\lambda\in[\lambda^*,\infty)$, $\underline{u}  =  \lambda^{r} \psi $ with $r = 1/(p+\beta-1)$,    is a lower solution of  $(P)_\lambda$.
\end{theorem}

\noindent {\bf Proof of Theorem \ref{teo 0.2}}  By $(f)_2\mbox{\rm  (i)-(ii)}$ there is $b>0$ such that
\begin{equation}\label{4.0}
f(s)>-\frac{b}{s^\beta}~\mbox{for}~ s>0.
\end{equation}

\nd  By lemma \ref{lem 3.1} there are both  a function $\phi \in C^{1, \alpha}(\overline{\Omega})$,  with $\alpha \in (0,1)$, such that
\begin{equation}\label{eq 4.1}
 \left\{\begin{array}{rllr}
-\Delta_{p} \phi&=& \frac{1}{\phi^\beta}~ \mbox{in}~ \Omega,\\
\phi&>&0~  \mbox{in}~  {\Omega,} \\
\phi&=&0~  \mbox{on}~ \partial {\Omega} ,
\end{array}
\right.
\end{equation}
\nd and  a constant $C_1>0$ such that $\phi \geq C_1  d~\mbox{in}~ \Omega$. Take  $\delta =a^{\frac{p-1}{\beta-1+p}}$ and $\gamma=2^\beta b\delta^{-\frac{\beta}{p-1}}$, where $a$ is given in $(f)_2 (i)$.
\vskip.3cm

\nd  By  corollary \ref{cor 3.2} there is a constant $\epsilon_0 > 0$ such that  for eaxh $\epsilon \in (0, \epsilon_0)$, the problem
\begin{equation}\label{eq 4.2}
\left\{\begin{array}{cclr}
-\Delta_p \psi & = & \delta\phi^{-\beta}\chi_{[d>\epsilon]}-\gamma\phi^{-\beta}\chi_{[d<\epsilon]} & \mbox{in}\quad \Omega,\\
\psi        & > & 0  & \mbox{in}\quad \Omega,\\
  \psi         & = & 0  & \mbox{on}\quad \partial\Omega,
\end{array}\right.
\end{equation}
\nd admits a solution $\psi \in C^{1, \alpha}(\overline{\Omega})$ satisfying
\begin{equation}\label{eq 4.3}
\psi \geq  ({{\delta^{{1}/{(p-1)}}}} / {2})\phi.
\end{equation}  
\nd Set
$\underline{u} =\lambda^r\psi~\mbox{where}~ r={1}/{(p+\beta-1)}~ \mbox{and}~\lambda >0$. Take $\lambda_* = [{2A}/(C_1 \epsilon\delta^{\frac{1}{p-1}})]^{\frac{1}{r}}$,
with  $\epsilon \in (0, \epsilon_0)$  and $A$   given by $(f)_2$. 
\vskip.3cm

\nd {\bf Claim}~~ $\underline{u}$ is a lower solution of $(P)_{\lambda}$ for $\lambda \geq \lambda_*$.
\vskip.2cm

\nd Indeed,  take $\xi\in\w$, $\xi\geq 0$. Using  (\ref{eq 4.2})    we have
\begin{equation}\label{eq 4.5}
\int_{\Omega}|\nabla \underline{u}|^{p-2}\nabla\underline{u}.\nabla \xi dx{=}\lambda^{r(p-1)}\delta\int_{\{d>\epsilon\}}\frac{\xi}{\phi^{\beta}}dx-\lambda^{r(p-1)}\gamma\int_{\{d<\epsilon\}}\frac{\xi}{\phi^{\beta}}dx.
\end{equation}

\nd We distinguish between two cases.
\vskip.1cm

\nd {\bf Case 1}~~ $d>\epsilon$
\vskip.1cm

\nd For each $\lambda\geq\lambda_*$ we have by using (\ref{eq 4.3}),
$$\displaystyle\underline{u}=\lambda^{r}\psi{\geq}\lambda^{r}\frac{\delta^{\frac{1}{p-1}}}{2}\phi \geq \lambda^{r}\frac{\delta^{\frac{1}{p-1}}}{2}C_{1}d >\lambda^{r}\frac{\delta^{\frac{1}{p-1}}}{2}C_{1}\epsilon>A.$$

\nd So $\underline{u}(x)>A$ for each $\lambda\geq\lambda_*$ with $d(x) > \epsilon$.  By (\ref{eq 4.1}) and (\ref{eq 4.2}),
$$-\Delta_p {\delta^{\frac{1}{p-1}}\phi }~{=}~ \frac{\delta}{\phi^\beta}~{\geq}~-\Delta_p\psi.$$
\nd It follows by the weak comparison principle that
\begin{equation}\label{eq 4.4}
\delta^{\frac{1}{p-1}}\phi\geq \psi \quad \mbox{in} \ \Omega.
\end{equation}
\nd Using $(f)_2 \mbox{\rm (i)}$ and  (\ref{eq 4.4}) we have,
\begin{eqnarray}\label{eq 4.6}
\displaystyle\lambda \int_{d>\epsilon}f(\underline{u})\xi dx &{\geq} & \displaystyle\lambda a \int_{d>\epsilon}\frac{\xi}{{\underline{u}}^{\beta}} dx
{=} \displaystyle\lambda^{1-r\beta}a\int_{d>\epsilon}\frac{\xi}{{\psi}^{\beta}}dx \nonumber\\
&{\geq} &\displaystyle\lambda^{r(p-1)}
\frac{a}{\delta^{\frac{\beta}{p-1}}}\int_{d>\epsilon}\frac{\xi}{\phi^{\beta}}dx= \displaystyle\lambda^{r(p-1)}\delta\int_{d>\epsilon}\frac{\xi}{\phi^{\beta}}dx.
\end{eqnarray}
\vskip.2cm

\nd {\bf Case 2}   $d<\epsilon$.
\vskip.1cm

\nd Using (\ref{4.0}) and (\ref{eq 4.3})  we have

\begin{eqnarray} \label{eq 4.7}
\displaystyle\lambda \int_{\{d<\epsilon\}}f(\underline{u})\xi dx &\geq & \displaystyle-\lambda b \int_{\{d<\epsilon\}}\frac{\xi}{{\underline{u}}^{\beta}}dx {=}\displaystyle-\lambda^{1-r\beta}b\int_{d<\epsilon}\frac{\xi}{{\psi}^{\beta}}dx\nonumber\\
 &{\geq}&\displaystyle-\lambda^{r(p-1)}b
\frac{2^{\beta}}{\delta^{\frac{\beta}{p-1}}}\int_{d<\epsilon}\frac{\xi}{\phi^{\beta}}dx=-\lambda^{r(p-1)}\gamma\int_{d<\epsilon}\frac{\xi}{\phi^{\beta}}dx.
 \end{eqnarray}
 \nd Using  (\ref{eq 4.6})-(\ref{eq 4.7}) we get
 $$
 \lambda\int_\Omega f(\underline{u})\xi dx +\int_\Omega h\xi dx \geq \int_{\Omega}|\nabla \underline{u}|^{p-2}\nabla\underline{u}.\nabla \xi dx,
 $$
\nd showing that  $\underline{u}=\lambda^{r}\psi$ is a lower solution of  $(P)_\lambda$ for each $\lambda \geq\lambda_*$, ending the proof of theorem
\ref{teo 0.2}.  $\hfill{\rule{2mm}{2mm}}$
\vskip.3cm

\nd Next, we show  existence of an upper solution.

\begin{theorem} \label{teo 0.3}
Assume $(f)_1-(f)_2$ and take $\Lambda>\lambda_*$ with $\lambda_*$ as in theorem $\ref{teo 0.2}$. Then for each  $\lambda\in[\lambda^*,\Lambda]$, $(P)_\lambda$ admits an upper solution $\overline{u}=\overline{u}_\lambda = M \phi$
  where $M > 0$ is a constant and $\phi$ is given by $(\ref{eq 4.1})$.
\end{theorem}

\noindent {\bf Proof of Theorem \ref{teo 0.3}} Choose  $\overline{\epsilon}>0$ such that
 \begin{equation}\label{eq 4.8}
 \displaystyle\Lambda\overline{\epsilon}||\phi||^{p-1+\beta}_{\infty}<\frac{1}{4}.
 \end{equation}
\nd By  $(f)_1$ and  $(f)_2$ there are $A_1>0$ and $C>0$ such that
\begin{equation}\label{eq 4.9}
|f(u)|\leq \overline{\epsilon}u^{p-1}~ \mbox{for}~ u>A_{1}
\end{equation}
\nd and
\begin{equation}\label{eq 4.10}
|f(u)|\leq \frac{C}{u^{\beta}}~\mbox{for}~ u \leq A_{1}.
\end{equation}
\nd Choose
\begin{equation}\label{eq 4.11}
M \geq\bigg\{\Lambda^{r}\delta^{\frac{1}{p-1}}, \ (4\Lambda C)^{\frac{1}{p+\beta-1}},\big(4||h||_\infty||\phi||_\infty^\beta\big)^{\frac{1}{p-1}}\bigg\}.
\end{equation}
\nd Using  (\ref{eq 4.8}) and (\ref{eq 4.11})  we get
\begin{equation}\label{eq 4.12}
\displaystyle\Lambda \overline{\epsilon}\big(M||\phi||_{\infty}\big)^{p+\beta-1}+\Lambda C\leq \frac{M^{p+\beta-1}}{4}+\frac{M^{p+\beta-1}}{4} =\frac{M^{p+\beta-1}}{2}.
\end{equation}
\nd Let  $\overline{u}=M\phi$. Using  (\ref{eq 4.9})-(\ref{eq 4.10}) and picking  $\lambda\leq\Lambda$ we have
\begin{eqnarray}\label{eq 4.13}
\lambda f(\overline{u})& \leq &\lambda|f(\overline{u})|\nonumber\\
& \leq & \lambda\bigg[\overline{\epsilon}\, \overline{u}^{p-1}\chi_{\{\overline{u}>A_1\}}+\frac{ C}{\overline{u}^\beta}\chi_{\{\overline{u}\leq A_1\}}\bigg]\nonumber\\
& \leq &  \lambda\bigg[\overline{\epsilon}\, \overline{u}^{p-1}\chi_{\{\overline{u}>A_1\}}+\overline{\epsilon}\, \overline{u}^{p-1}\chi_{\{\overline{u}\leq A_1\}}+\frac{ C}{\overline{u}^\beta}\chi_{\{\overline{u}\leq A_1\}}+\frac{ C}{\overline{u}^\beta}\chi_{\{\overline{u}> A_1\}}\bigg]\nonumber\\
& = & \lambda\bigg[\overline{\epsilon}\, \overline{u}^{p-1}+\frac{ C}{\overline{u}^\beta}\bigg].\nonumber\\
\end{eqnarray}
\nd Thus
\begin{eqnarray}\label{eq 4.14}
\lambda f(M\phi) & \leq &  \lambda\bigg[\frac{\overline{\epsilon}(M||\phi||_{\infty})^{p+\beta-1}+C}{[M \phi]^{\beta}}\bigg]\nonumber\\
& \leq & \Lambda\frac{\overline{\epsilon}(M||\phi||_{\infty})^{p+\beta-1}}{[M \phi]^{\beta}}+\Lambda\frac{C}{[M \phi]^{\beta}}.
\end{eqnarray}
\nd Replacing (\ref{eq 4.11}) and (\ref{eq 4.12}) in (\ref{eq 4.14}),
$$
\lambda f(M\phi)  \leq   \frac{M^{p+\beta-1}}{2[M {\phi}]^{\beta}}  =   \frac{M^{p-1}}{2\phi^{\beta}}.
$$
\nd It follows from  (\ref{eq 4.11}) that
$$
h  \leq  ||h||_\infty\\
      \frac{M^{p-1}}{2||\phi||_\infty^{\beta}}\\
      \leq  \frac{M^{p-1}}{2\phi^\beta}.
$$
\nd Thus
$$\lambda f(\overline{u})+h \leq \frac{M^{p-1}}{\phi^\beta}.$$
\vskip.1cm

\nd Taking $\eta\in\w$ with $\eta\geq0$ we have by using (\ref{eq 4.1}),
\begin{eqnarray*}
\lambda \int_{\Omega}f(\overline{u})\eta dx + \int_\Omega h\eta dx&\leq& M^{p-1}\int_{\Omega}\frac{\eta}{\phi^{\beta}}dx\\
&{=}& M^{p-1}\int_{\Omega}|\nabla \phi|^{p-2}\nabla \phi.\nabla \eta dx\\
&=& \int_{\Omega}|\nabla (M\phi)|^{p-2}\nabla (M\phi).\nabla \eta dx\\
&=& \int_{\Omega}|\nabla \overline{u}|^{p-2}\nabla \overline{u}.\nabla \eta dx,
\end{eqnarray*}
\nd showing that  $\overline{u}=M\phi$ is an upper solution of $(P)_\lambda$ for $\lambda\in[\lambda_*,\Lambda]$.  $\hfill{\rule{2mm}{2mm}}$

\section{Proofs of the Main Results}

\nd At first we introduce some notations, remarks and lemmas.  Take $\Lambda>\lambda_*$ and set   $I_{\Lambda} :=[\lambda_*,\Lambda]$.  For each $\lambda \in  I_{\Lambda} $. By theorem \ref{teo 0.2}
$$
\underline{u}=\underline{u}_\lambda = \lambda^{r}\psi
$$
\nd is a lower solution of $(P)_\lambda$. Pick  $M = M_{\Lambda} \geq \Lambda^{r}\delta^{\frac{1}{p-1}}$.  By theorem \ref{teo 0.3}
$$
\overline{u}=\overline{u}_\lambda = M_{\Lambda} \phi
$$
\nd is an upper solution of $(P)_\lambda$. It follows by (\ref{eq 4.4}) that
\begin{equation}\label{sub-menor-super}
\underline{u} = \lambda^{r}\psi \leq   \Lambda^{r}\delta^{\frac{1}{p-1}}\phi \leq  M\phi  = \overline{u}.
\end{equation}

\nd The convex, closed subset of $I_{\Lambda} \times C(\overline{\Omega})$, defined by

$$
\mathcal{G}_{\Lambda} := \big\{(\lambda,u) \in I_{\Lambda} \times C(\overline{\Omega})~|~ \lambda \in I_{\Lambda},~\underline{u}\leq u\leq \overline{u}\ \mbox{and}\ u=0 \ \mbox{on}\ \partial\Omega\big\}
$$
\nd will play a key role in this work.
\vskip.2cm

\nd For each  $u\in C(\overline{\Omega})$ define
\begin{equation}\label{eq 4.15}
{f}_{\Lambda}(u)=\chi_{S_1}f(\underline{u})+\chi_{S_2}f(u)+\chi_{S_3}f(\overline{u}),~ x\in\Omega,
\end{equation}
\nd where
$$
\begin{array}{lcl}
S_1:=\big\{x\in \Omega~|~u(x)<\underline{u}(x)\big\},\\

S_2:=\big\{x\in \Omega~|~ \underline{u}(x)\leq u(x)\leq \overline{u}(x)\big\},\\

S_3:=\big\{x\in \Omega~|~ \overline{u}(x)<u(x)\big\},
\end{array}
$$
\nd and  $\chi_{S_i}$  is the characterictic function of  $S_i$.
\vskip.1cm

\begin{lemma}\label{lem 5.1}
For each $u \in C(\overline{\Omega})$, ${f}_{\Lambda}(u)\in L^\infty_{loc}(\Omega)$ and there are  $C>0$, $\beta\in(0,1)$ such that
\begin{equation}\label{DES f LAMBDA}
|{f}_{\Lambda}(u)(x)|\leq\frac{C}{d(x)^\beta},~~ x \in \Omega.
\end{equation}
\end{lemma}

\nd {\bf Proof }  Indeed,  let $\mathcal{K}\subset\Omega$ be a compact subset. Then both    $\underline{u}$ and  $\overline{u}$ achieve a positive  maximum and a positive minimum  on  $\mathcal{K}$. Since   $f$ is continuous in $(0,\infty)$ then ${f}_{\Lambda}(u)\in L^\infty_{loc}(\Omega)$.
\vskip.2cm

\nd {\bf Verification of (\ref{DES f LAMBDA}):} Since $\displaystyle \Omega = \cup_{i=1}^{3} S_i$ it is enough to show that
$$
|f(u(x))| \leq \frac{C}{d(x)^{\beta}},~~ x \in S_i,~  i = 1,2,3.
$$
\nd At first, by $(f)_{2}\mbox{\rm (ii)}$  there are $C, \delta > 0$ such that
$$
|f(s)| \leq \frac{C}{s^{\beta}},~~ 0 < s < \delta.
$$
\nd Let
$$
\Omega_{\delta} = \{x \in \Omega~|~d(x) < \delta \}.
$$
\nd   Recalling that  $\underline{u} \in C^{1}(\overline{\Omega})$,  let
$$
D = \max_{\overline{\Omega}}  d(x),~~     \nu_\delta := \min_{\overline{\Omega_{\delta}^{c}}} d(x),~~\nu^\delta := \max_{\overline{\Omega_{\delta}^{c}}} d(x),
$$
\nd  and notice that both $0 < \nu_\delta \leq \nu^\delta \leq D  < \infty$ and $f([\nu_\delta, \nu^\delta])$
is compact.
\vskip.2cm

 \nd On the other hand, applying theorems \ref{teo 0.2}, \ref{teo 0.3},  lemmas \ref{lem 3.2}, \ref{lem 3.1} and inequality (\ref{eq 4.3}) we infer that
$$
0 < \lambda_*^{r} \psi = \underline{u}~  \leq~ \overline{u} =  M \phi~\mbox{in}~\Omega
$$
\nd and
$$
\frac{1}{\underline{u}^{\beta}} ,~  \frac{1}{\overline{u}^{\beta}}, \leq \frac{1}{ (\lambda_*^{r} \psi(x))^{\beta}} \leq \frac{C}{d(x)^{\beta}},~ x \in \Omega_{\delta}.
$$
\nd  To finish to proof,  we distinguish among  three cases:
\vskip.3cm

\nd  { {\bf (i)}  $x \in S_1$:}  in this case,
$$
f_{\Lambda}(u(x)) = f(\underline{u}(x)).
$$
\nd If $x \in S_1 \cap \Omega_{\delta}$ we infer that
$$
|f_{\Lambda}(u(x))|  \leq \frac{C}{\underline{u}(x)^{\beta}} \leq \frac{C}{d(x)^{\beta}}.
$$
\nd If  $x \in S_1 \cap \Omega_{\delta}^c$. Pick  positive numbers $ d_i,~i=1,2$ such that
$$
 d_1 \leq  \underline{u}(x) \leq d_2,~x \in \Omega_{\delta}^c.
$$
\nd Hence
$$
|f_{\Lambda}(u(x))| \leq \frac{C}{d(x)^{\beta}},~x \in \Omega.
$$
\nd {\bf (ii)  $x \in S_2$:}  in this case,
$$
0 < \lambda_*^{r} \psi \leq u  \leq  M \phi.
$$
\nd and as a consequence,
$$
|f(u(x))| \leq \frac{C}{u(x)^{\beta}},~~x \in \Omega_{\delta}.
$$
\nd Hence, there is a positive constant $\widetilde{C}$ such that
$$
|f(u(x))| \leq \widetilde{C},~~x \in \overline{\Omega_{\delta}^{c}}.
$$
\nd Thus
$$
|f(u(x))| \leq \left\{
\begin{array}{rclr}
\widetilde{C}~ &\mbox{if}& x \in \overline{\Omega_{\delta}^{c}},\\
\frac{C}{d(x)^{\beta}}~ &\mbox{if}& x \in  \Omega_{\delta}.
\end{array}
\right.
$$
\nd On the other hand,
$$
\frac{1}{D^{\beta}}  \leq  \frac{1}{d(x)^{\beta}},~ x \in \overline{\Omega_{\delta}^{c}},
$$
\nd and therefore there is a constant $C > 0$ such that
$$
|f(u(x))| \leq \left\{
\begin{array}{rclr}
\frac{C}{D^\beta}~ &\mbox{if}& x \in \overline{\Omega_{\delta}^{c}},\\
\\
\frac{C}{d(x)^{\beta}}~ &\mbox{if}& x \in  \Omega_{\delta}.
\end{array}
\right.
$$
\nd Therefore,
$$
|f(u(x))| \leq \frac{C}{d(x)^{\beta}},~~x \in S_2,~u \in \mathcal{G}_{\Lambda}.
$$
\nd {\bf Case $x \in S_3$:}  in this case
$$
f_{\Lambda}(u(x)) = f(\overline{u}(x)).
$$
\nd If $x \in S_3 \cap \Omega_{\delta}$ we infer that
$$
|f_{\Lambda}(u(x))|  \leq \frac{C}{\overline{u}(x)^{\beta}} \leq \frac{C}{d(x)^{\beta}}.
$$
\nd If  $x \in S_3 \cap \Omega_{\delta}^c$. Pick  positive numbers $ d_i,~i=1,2$ such that
$$
  d_1 \leq \overline{u}(x)) \leq d_2,~ x \in  \Omega_{\delta}^c.
$$
\nd Hence
$$
|f_{\Lambda}(u(x))| \leq \frac{C}{d(x)^{\beta}},~x \in \Omega.
$$
\nd This ends the proof of  lemma \ref{lem 5.1}. $\hfill{\rule{2mm}{2mm}}$

\begin{remark}
  By lemmas $\ref{lem 3.3}, \ref{lem 5.1}$ and remark  $(\ref{SOL OPER})$, for each $v \in C(\overline{\Omega)}$ and $\lambda \in I_{\Lambda}$,
\begin{equation}\label{EST f LAMBDA}
 (\lambda {f}_{\Lambda}(v)+h)  \in L^\infty_{loc}(\Omega)~\mbox{and}~|(\lambda {f}_{\Lambda}(v)+h)|\leq \frac{C_{\Lambda}}{d^\beta(x)} ~ \mbox{ in}~\Omega
\end{equation}
\nd where $C_{\Lambda} >0 $ is a constant independent of $v$ and $\beta\in(0,1)$.  So for each  $v$,
\begin{equation}\label{eq 4.15}
\left\{
\begin{array}{rllr}
-\Delta_{p} u&=& \lambda {f}_{\Lambda}(v)+h~ \mbox{in} \quad \Omega,\\

u&=&0 ~ \mbox{on} ~ \partial {\Omega}
\end{array}
\right.
\end{equation}
\nd   admits an only solution   $u = S(\lambda {f}_{\Lambda} (v)+h) ) \in \w\cap C^{1,\alpha}(\overline{\Omega})$.
\end{remark}
\nd Set
$$
{F}_{\Lambda}(u)(x)={f}_{\Lambda}(u(x)),~  u\in C(\overline{\Omega}).
$$
\nd and consider the operator
$$
T:I_{\Lambda}\times C(\overline{\Omega}) \to  \w \cap C^{1,\alpha}(\overline{\Omega}),
$$

$$
T(\lambda,u) =
S(\lambda {F}_{\Lambda} (u)+h) )~\mbox{if}~  \lambda_* \leq \lambda \leq \Lambda, ~ u \in C(\overline{\Omega}).
$$

\nd  Notice that if $(\lambda, u) \in I_{\Lambda} \times C(\overline{\Omega})$ satisfies $u = T(\lambda,u)$ then $u$ is a solution of
\begin{equation}
\left\{
\begin{array}{rllr}
-\Delta_{p} u&=& \lambda {f}_{\Lambda}(u)+h~ \mbox{in} \quad \Omega,\\

u&=&0 ~ \mbox{on} ~ \partial {\Omega}
\end{array}
\right.
\end{equation}
\begin{lemma}\label{lem 4.2}
If $(\lambda,u) \in I_{\Lambda} \times C(\overline{\Omega})$ and  $u = T(\lambda,u)$    then  $(\lambda,u) \in\mathcal{G}_{\Lambda}$.
 \end{lemma}
 \prova
 Indeed, let  $(\lambda, u) \in I_{\Lambda} \times C(\overline{\Omega})$  such that $T(\lambda,u) = u$. Then
$$
\int_\Omega|\nabla u|^{p-2}\nabla u.\nabla\varphi dx=\lambda\int_\Omega {f}_{\Lambda}(u)\varphi dx +\int_\Omega h\varphi dx,~ \varphi\in\w.
$$
\nd We claim that  $u\geq \underline{u}$. Assume on the contrary, that $\varphi :=(\underline{u}-u)^{+}\not\equiv 0$. Then
$$\begin{array}{rcl}
\displaystyle \int_{\Omega}|\nabla u|^{p-2}\nabla u.\nabla\varphi dx &= &\displaystyle\int_{u<\underline{u}}|\nabla u|^{p-2}\nabla u.\nabla\varphi dx\\
  &=&\displaystyle\lambda\int_{u<\underline{u}}\ {f}_{\Lambda}(u).\varphi dx + \int_{u<\underline{u}} h\varphi dx \\
   &=&\displaystyle\lambda\int_{u<\underline{u}}\ f(\underline{u}).\varphi dx + \int_{u<\underline{u}} h\varphi dx\\
&  \geq& \displaystyle\int_{u<\underline{u}}|\nabla \underline{u}|^{p-2}\nabla \underline{u}.\nabla\varphi dx\\
& = &\displaystyle \int_{\Omega}|\nabla \underline{u}|^{p-2}\nabla \underline{u}.\nabla\varphi dx.
\end{array}
$$
\nd Hence
$$
\displaystyle\int_{\Omega}\big [|\nabla u|^{p-2}\nabla u-|\nabla \underline{u}|^{p-2}\nabla\underline{u} \big] \cdot  \nabla(u-\underline{u})  dx \leq 0.
$$
\nd It follows by  lemma \ref{lem A.1}  that
$\displaystyle\int_{\Omega}|\nabla \big(u-\underline{u} \big )|^{p} dx \leq 0$, contradicting   $\varphi \not\equiv 0$. Thus,  $(\underline{u}-u)^{+}=0$, that is,  $\underline{u}-u\leq0$,  and so $\underline{u}\leq T(\lambda,u)$.
\vskip.2cm

\nd We claim that  $\overline{u}\geq u$. Assume on the contrary that  $\varphi :=(u - \overline{u})^{+}\not\equiv 0$. We have
$$
\begin{array}{rcl}
\displaystyle \int_{\Omega}|\nabla u|^{p-2}\nabla u.\nabla\varphi dx &= &\displaystyle\int_{\overline{u}<u}|\nabla u|^{p-2}\nabla u.\nabla\varphi dx\\
  &=&\displaystyle\lambda\int_{\overline{u}<u}\ {f}_{\Lambda}(u).\varphi dx + \int_{\overline{u}<u}h\varphi dx \\
   &=&\displaystyle\lambda\int_{\overline{u}<u}\ f(\overline{u}).\varphi dx + \int_{\overline{u}<u}h\varphi dx\\
&  \leq& \displaystyle\int_{\overline{u}<u}|\nabla \overline{u}|^{p-2}\nabla \overline{u}.\nabla\varphi dx\\
& = &\displaystyle \int_{\Omega}|\nabla \overline{u}|^{p-2}\nabla \overline{u}.\nabla\varphi dx,
\end{array}
$$
\nd Therefore,
$$
\displaystyle\int_{\Omega}\big [|\nabla u|^{p-2}\nabla u-|\nabla \overline{u}|^{p-2}\nabla\overline{u} \big]    \cdot  \nabla(u-\overline{u})  dx \leq 0.
$$
\nd contradicting $\varphi  \not\equiv 0$. So  $(u-\overline{u})^{+}=0$ so that $u-\overline{u}\leq0$, which gives  $\overline{u}\geq T(\lambda,u)$.
\vskip.1cm

\nd As a consequence of the arguments above $u\in\mathcal{G}_{\Lambda}$, showing  lemma \ref{lem 4.2}.  $\hfill{\rule{2mm}{2mm}}$
\vskip.2cm

\begin{remark}\label {REMARK X1} By the definitions of $f_{\Lambda}$ and $\mathcal{G}_{\Lambda}$, for each $(\lambda,u) \in \mathcal{G}_{\Lambda}$
\begin{equation}\label{eq 4.15}
{f}_{\Lambda}(u)=f(u),~~ x\in\Omega.
\end{equation}
\end{remark}
\begin{remark}\label{REMARK X2}   By remark $\ref{SOL OPER}$,  there is $R_{\Lambda} > 0$ such that
$\mathcal{G}_{\Lambda} \subset B(0,R_{\Lambda}) \subset  C(\overline{\Omega})$  and
$$
T \Big (I_{\Lambda} \times \overline{B(0,R_{\Lambda})}\Big) \subseteq  {B(0,R_{\Lambda})} .
$$
\nd Notice that, by $(\ref{eq 4.15})$ and lemma $\ref{lem 4.2}$,   if $(\lambda,u) \in I_{\Lambda} \times C(\overline{\Omega})$ satisfies $u = T(\lambda,u)$ then $(\lambda,u)$ is a solution of  ${(P)_\lambda}$. \nd By remark $\ref{REMARK X1}$,  to solve $(P)_{\lambda}$ it suffices to  look for fixed points of $T$.
\end{remark}

\begin{lemma}\label{lem 4.1}~ $T: I_{\Lambda}  \times \overline{B(0,R_{\Lambda})} \to  \overline{B(0,R_{\Lambda})}$ is continuous and compact.
\end{lemma}

\nd \prova Let  $\{(\lambda_n, u_{n})\} \subseteq I_{\Lambda}  \times \overline{B(0,R_{\Lambda})} $ be a  sequence such that
$$
\lambda_n \to \lambda~~\mbox{and}~~u_n \stackrel{ C(\overline{\Omega})  }  \rightarrow  u.
$$
\nd Set
$$
v_n = T(\lambda_n, u_{n})~\mbox{and}~ v = T(\lambda, u)
$$
\nd so that
$$
v_n = S(\lambda_n {F_{\Lambda}}(u_n)+h)~\mbox{ and}~
v = S(\lambda {F_{\Lambda}}(u)+h).
$$
\nd It follows that
$$
\begin{array}{lcl}
 \displaystyle \int_\Omega\big(|\nabla v_n|^{p-2}\nabla v_n-|\nabla v|^{p-2} \nabla v\big).\nabla(v_n-v)dx &=&\displaystyle\lambda_n\int_\Omega\big({f_{\Lambda}}(u_n)- {f_{\Lambda}}(u)\big)(v_n-v) dx\\
 & \leq & \displaystyle C \int_\Omega| {f_{\Lambda}}(u_n)- {f_{\Lambda}}(u)|dx.
 \end{array}
 $$
\nd Since
$$
|f_{\Lambda}(u_n)-f_{\Lambda}(u)|\leq \frac{C}{d(x)^\beta} \in L^1(\Omega)~\mbox{and}~
f_{\Lambda}(u_n(x))\to f_{\Lambda}(u(x))~\mbox{a.e.}~x \in \Omega,
$$
\nd it follows by Lebesgue's  Theorem that
$$
\int_\Omega|f_{\Lambda}(u_n)-f_{\Lambda}(u)|dx\to 0.
$$
\nd Therefore $v_n \to v$ in $W_0^{1,p}(\Omega)$.
\vskip.2cm

\nd On the other hand, since  $u_n \stackrel{ C(\overline{\Omega})  }  \longrightarrow  u$, by the proof of lemma \ref{lem 5.1},
$$
 (\lambda_n {f}_{\Lambda}(u_n)+h)  \in L^\infty_{loc}(\Omega)~\mbox{and}~|(\lambda_n {f}_{\Lambda}(u_n)+h)|\leq \frac{C_{\Lambda}}{d^\beta(x)} ~ \mbox{ in}~\Omega.
$$
\nd By lemma \ref{lem 3.3} there is a constant $M > 0$ such that
$$
||v_n||_{C^{1,\alpha}(\overline{\Omega})} \leq M
$$
\nd so that $v_n \stackrel{ C(\overline{\Omega})  }  \rightarrow v$. This shows that $T: I_{\Lambda}  \times \overline{B(0,R_{\Lambda})} \to  \overline{B(0,R_{\Lambda})}$ is  continuous.
\vskip.1cm

\nd The  compactness of $T$ follows from the arguments in the five lines above.  $\hfill{\rule{2mm}{2mm}}$

\subsection  {\bf \large Proof of Theorem \ref{teo 0.0}}

\nd  Some notations and technical  results are needed.  At first, we recall  the Leray-Schauder Continuation Theorem (see \cite{KD},\cite{CG}).
\vskip.2cm

\begin{theorem}\label{teo A1}
\nd Let  $D$ be an open bounded subset of the Banach space $X$. Let $a,b \in {\bf R}$ with $a < b$
 and assume that  $T:[a,b] \times \overline{D} \rightarrow X$ is compact and continuous. Consider
$\Phi:[a,b] \times \overline{D} \rightarrow X$ defined by $\Phi(t,u)=u-T(t,u)$.
\nd Assume that
$$
\mbox{\rm(i)}~~  \displaystyle\Phi(t,u)\neq 0,~  t\in [a.b],~ u \in \partial D,~~
\mbox{\rm(ii)}~~  \displaystyle  {\rm deg} \big(\Phi(t,.),D,0 \big)\neq 0~ \mbox{for some}~  t \in [a,b].
$$
\nd and set
$$
\mathcal{S}_{a,b} = \{(t,u) \in [a,b] \times \overline{D} ~|~ \Phi(t,u) = 0\}.
$$
\nd Then, there is a connected compact subset $\Sigma_{a,b}$ of $\mathcal{S}_{a,b}$ such that

$$
\Sigma_{a,b}\cap (\{a\}\times D)\neq \emptyset
$$
\nd and
$$
\Sigma_{a,b}\cap (\{b\}\times D)\neq\emptyset.
$$

\end{theorem}

\nd The Leray-Schauder Theorem above will be applied to the operator  $T$ in the settings of {\it Section 5}. Remember that $T$ continuous, compact and $T \big ( I_{\Lambda} \times \overline{B(0,R_{\Lambda})} \big ) \subset {B(0,R_{\Lambda})}$. Consider  $\displaystyle \Phi : I_{\Lambda} \times \overline{B(0,R)}\!\longrightarrow \! \overline{B(0,R)})$ defined by
$$
\Phi(\lambda,u)=u-T(\lambda,u).
$$
\begin{lemma}\label{lem 4.1a} $\Phi$ satisfies:
\end{lemma}
\begin{itemize}
\item[(i)] $\Phi(\lambda,u)\neq 0$~ $ (\lambda,u)\in I_{\Lambda} \times \partial B(0,R_{\Lambda})$,
\item[(ii)] $\displaystyle  {\rm deg} (\Phi(\lambda,.),B(0,R_{\Lambda}),0)\neq 0$~\mbox{for each}~  $\lambda \in I_{\Lambda}$,
\end{itemize}\vspace{0,2 cm}
\nd \prova The verification of (i) is straightforward since  $T \big ( I_{\Lambda} \times \overline{B(0,R_{\Lambda})} \big ) \subset {B(0,R_{\Lambda})}$.
\vskip.2cm

\nd To prove (ii) , set $R = R_{\Lambda}$,  take  $\lambda \in I_{\Lambda}$ and consider the homotopy
$$
\Psi_\lambda(t,u)=u-tT(\lambda,u),~ (t,u)  \in [0,1] \times \overline{B(0,R)}.
$$
\nd It follows that  $0 \notin \Psi_\lambda(I\times\partial B(0, R))$. By the invariance under homotopy property of the Leray-Schauder degree
$$
{\rm deg} (\Psi_\lambda(t,.),B(0, R),0)= {\deg} (\Psi_\lambda(0,.),B(0,R),0)=1,~  t\in[0,1].
$$
\nd Setting
$$
\Phi(\lambda,u)=u-T(\lambda,u),~ (\lambda,u)  \in I_{\Lambda} \times \overline{B(0,R)},
$$
\nd we also have
$$
{\rm deg} (\Phi(\lambda,.),B(0,R),0)=1, \quad \lambda\in I_{\Lambda}.
$$
\nd Set
 $$
 \mathcal{S}_\Lambda =\big\{(\lambda,u)\in I_{\Lambda} \times \overline{B(0,R)} ~|~ \Phi(\lambda,u)=0\big\} \subset \mathcal{G}_{\Lambda}.
$$
\nd By the Leray-Schauder Continuation Theorem, there is a connected component $ \Sigma_{\Lambda} \subset \mathcal{S}_\Lambda$ such that
$$
\Sigma_{\Lambda} \cap (\{\lambda_*\}\times \overline{B(0,R)} ) \neq\emptyset
$$
\nd and
$$
\Sigma_{\Lambda} \cap (\{\Lambda\}\times \overline{B(0,R)} ) \neq\emptyset.
$$
\nd We point out that $\mathcal{S}_\Lambda$ is the solution set of the auxiliary problem
$$
\left\{
\begin{array}{rllr}
-\Delta_{p} u&=& \lambda {f}_{\Lambda}(u)+h~ \mbox{in} \quad \Omega,\\

u&=&0~~ \mbox{on}~~ \partial {\Omega}
\end{array}
\right.
$$
\nd and since $ \Sigma_{\Lambda} \subset \mathcal{S}_\Lambda \subset \mathcal{G}_{\Lambda}$ it  follows using the definition of $f_{\Lambda}$ that
$$
\left\{
\begin{array}{rllr}
-\Delta_{p} u&=& \lambda {f}(u)+h~ \mbox{in} \quad \Omega,\\

u&=&0~~ \mbox{on}~~ \partial {\Omega}
\end{array}
\right.
$$
\nd for $(\lambda,u) \in \Sigma_{\Lambda} $, showing that $ \Sigma_{\Lambda} \subset \mathcal{ S}$.
\vskip.1cm

\nd This ends the proof of theorem \ref{teo 0.0}.  $\hfill{\rule{2mm}{2mm}}$

\subsection {\bf \large Proof of Theorem \ref{teo 0.1}}

\nd We shall employ topological arguments to construct a suitable connected component of  of the solution set $ \mathcal{ S}$ of $(P)_\lambda$. To this aim some notations are needed.
\vskip.3cm

\nd Let $M = (M,d)$ be a metric space and denote by  $\{ \Sigma_n \}$ be a sequence of connected components of $M$. The {\it upper limit} of $\{ \Sigma_n \}$ is defined by
$$
\displaystyle \overline{\lim}~ \Sigma_n  = \{u \in M~|~ \mbox{there is }~ (u_{n_i}) \subseteq \cup \Sigma_{n}~\mbox{with}~ u_{n_i} \in \Sigma_{n_i}~\mbox{and}~u_{n_i} \to u    \}.
$$
\begin{remark}
$\overline{\lim}~ \Sigma_n$ is a closed subset of $M$.
\end{remark}

\nd We shall apply  theorem 2.1 in Sun \& Song \cite{sun-song}, stated below for the reader's convenience.

\begin{theorem}\label{Song}
Let $M$ be a metric space and $\{\alpha_n\}, \{\beta_n \} \in {\bf R}$ be sequences satisfying
$$
 \cdots <  \alpha_n < \cdots < \alpha_1 < \beta_1 < \cdots < \beta_n <  \cdots
$$
\nd with
$$
\alpha_n \to -\infty~\mbox{and}~ \beta_n \to \infty.
$$
\nd Assume that $\{\Sigma_n^* \}$ is  a sequence of connected subsets of ${\bf R} \times M$ satisfying,
$$
\begin{array}{lcl}
\mbox{\rm (i)}~~~\Sigma_n^* \cap (\{\alpha_n\}\times M)\neq \emptyset, \\
\\
\mbox{\rm (ii)}~~\Sigma_n^* \cap (\{\beta_n \}\times M)\neq \emptyset,
\end{array}
$$
\nd for each  $n$. For each $\alpha, \beta \in (-\infty, \infty)$ with $\alpha < \beta$,
$$
\mbox{\rm (iii)}~~~\big(\cup \Sigma_{n}^* \big) \cap ([\alpha,\beta] \times M)~\mbox{is a relatively compact subset of}~ {\bf R} \times M.
$$
\nd Then there is a connected component $\Sigma^*$ of~ $\overline{\lim}~ \Sigma_n^*$ such that
$$
\Sigma^* \cap (\{\lambda \}\times M)\neq \emptyset~\mbox{for each}~ \lambda \in (\lambda_*, \infty).
$$
\end{theorem}
\vskip.2cm

\nd \prova {\bf of Theorem \ref{teo 0.1} (finished)}~ Consider $\Lambda$ as introduced in  Section 5 and take a sequence $\{\Lambda_n \}$ such that $\lambda_* <  \Lambda_1 < \Lambda_2 < \cdots$ with $\Lambda_n \to \infty$. Set $\beta_n = \Lambda_n$ and take a sequence $\{\alpha_n \} \subset {\bf R}$ such that $\alpha_n \to -\infty$ and $\cdots <  \alpha_n < \cdots < \alpha_1 < \lambda_* $.
\vskip.2cm

\nd Following the notations of Section 5 consider the sequence of intervals
$I_n = [\lambda_*, \Lambda_n]$. Set $M = C(\overline{\Omega})$  and let
$$
\mathcal{G}_{\Lambda_n} := \big\{(\lambda,u) \in I_{n} \times \overline{B}_{R_{n}}~|~ \underline{u}\leq u\leq \overline{u},~ u=0 \ \mbox{on}\ \partial\Omega\big\},
$$
\nd  where  $R_n = R_{\Lambda_n}$. Consider the sequence of compact operators
$$
T_n : [\lambda_*,\Lambda_n] \times {\overline{B}}_{R_{n}}  \to  {\overline{B}}_{R_{n}}
$$
\nd  defined by
$$
T_n(\lambda,u) =
S(\lambda {F}_{\Lambda_n} (u)+h) )~\mbox{if}~  \lambda_* \leq \lambda \leq \Lambda_n, ~ u \in {\overline{B}}_{R_{n}}.
$$
\nd Next consider the  extension of ${T}_n$, namely $\widetilde{T}_n : {\bf R} \times {\overline{B}}_{R_{n}} \to {\overline{B}}_{R_{n}}$ defined by
$$
\widetilde{T}_n(\lambda,u) = \left\{\begin{array}{lllll}
T_n(\lambda_*,u)~~ \mbox{if}~~ \lambda \leq \lambda_*,\\

T_n(\lambda,u)~~ \lambda_* \leq \lambda \leq \Lambda_n,\\

T_n( \Lambda_n,u)~~ \mbox{if}~~ \lambda \geq \Lambda_n.
\end{array}\right.
$$
\nd Notice that $\widetilde{T}_n$ is continuous, compact.
\vskip.2cm

\nd Applying  theorem \ref{teo A1}  to $\widetilde{T}_n : [\alpha_n, \beta_n] \times {\overline{B}}_{R_{n}} \to {\overline{B}}_{R_{n}}$ we get a compact connected  component $\Sigma_{n}^*$ of
$$
\mathcal{S}_n =\big\{(\lambda,u)\in [\alpha_n,\beta_n] \times {\overline{B}}_{R_{n}} ~|~ \Phi_n(\lambda,u)=0\big\}, $$
\nd where 
$$
\Phi_n(\lambda,u)=u-\widetilde{T}_n(\lambda,u).
$$
\nd Notice that $\Sigma_{n} ^*$ is also a connected subset of $ {\bf R} \times M$. By theorem \ref{Song} there is 
 a connected component $\Sigma^*$ of~ $\overline{\lim}~ \Sigma_n^*$ such that
$$
\Sigma^* \cap (\{\lambda \}\times M)\neq \emptyset~\mbox{for each}~ \lambda \in {\bf R}.
$$
\nd Set $\Sigma =( [\lambda_*, \infty) \times M ) \cap \Sigma^*$. Then $\Sigma \subset {\bf R} \times M$ is connected and  
$$
\Sigma \cap (\{\lambda \}\times M)\neq \emptyset,~~ \lambda_* \leq \lambda < \infty.
$$
\nd We claim that  $\Sigma   \subset  \mathcal{ S}$. Indeed, at first notice that 
\begin{equation}\label{extension}
{\widetilde{T}_{n+1}} {\Big |_ { \big([\lambda_*,\Lambda_n] \times {\overline{B}}_{R_{n}} \big) }  }  = {\widetilde{T}_{n}} {\Big |_ { \big([\lambda_*,\Lambda_n] \times {\overline{B}}_{R_{n}} \big) }  }  = {T}_n.
\end{equation}
\nd If $(\lambda,u) \in \Sigma$ with $\lambda > \lambda_*$,  there is a sequence $(\lambda_{n_i}, u_{n_i}) \in \cup \Sigma_{n}^*$ with $(\lambda_{n_i}, u_{n_i})  \in \Sigma_{n_i}^*$ such that $\lambda_{n_i} \to \lambda~\mbox{and}~u_{n_i} \to u$. Then  $u \in B_{R_N}$ for some integer $N > 1$. 
\vskip.2cm

\nd We can assume that $(\lambda_{n_i}, u_{n_i}) \in [\lambda_*, \Lambda_N] \times B_{R_N}$. On the other hand, by (\ref{extension}),
$$
u_{n_i} = T_{n_i}(\lambda_{n_i}, u_{n_i}) = T_{N}(\lambda_{n_i}, u_{n_i}).
$$
\nd Passing to the limit we get
$$
u = T_{N}(\lambda, u)
$$
\nd which shows that $(\lambda,u) \in \Sigma_{{N}}$ and so
$$
(\lambda, u) \in  \mathcal{ S}:=\big{\{} (\lambda,u)\in(0,\infty)\times C(\overline{\Omega})\ \big{|}\  \mbox{u is a solution of}~ (P)_\lambda \  \big{\}}.
 $$
\nd This ends the proof of theorem \ref{teo 0.1}.   $\hfill{\rule{2mm}{2mm}}$

\section{Appendix}

\nd In this section we present proofs of lemma \ref{lem 3.3}, corollary \ref{cor 3.2} and   recall some results referred to in the paper. We begin with   the Browder-Minty Theorem, (cf. Deimling \cite{KD}). Let $X$ be a real  reflexive  Banaxh space with dual space $X^*$.
A map $F:X\to X^*$ is monotone if
$$
\langle Fx-Fy,x-y\rangle\geq 0,~x, y \in X,
$$
$F$ is hemicontinuous if
 $$
F(x+ty) \stackrel{*} \rightharpoonup Fx ~ \mbox{as}~ t\to 0,
$$
\nd and $F$ is coercive if
 $$
\frac{\langle  Fx,x\rangle}{|x|}\to \infty \ \mbox{as}\ |x|\to\infty.
$$

\begin{theorem} \label{teo C.3}
Let $X$ be a real reflexive Banach space and let $F:X\to X^*$ ve  a monotone, hemicontinous and coercive operator. Then $F(X) = X^*$. Moreover, if $F$ is strictly monotone then it is a homeomorphism.
\end{theorem}

\nd The inequality below, (cf \cite{S}, \cite{P}), is very useful when dealing with the p-Laplacian.
\begin{lemma}\label{lem A.1}
\nd Let $p>1$. Then there is a constant $C_p>0$ such that
\begin{equation}\label{simon}
\big(|x|^{p-2}x-|y|^{p-2}y,x-y \big)\geq \left\{\begin{array}{rl}
C_p~|x-y|^p & \mbox{\rm if}\quad p\geq 2,\\

C_p~\frac{|x-y|^p}{(1+|x|+|y|)^{2-p}} & \mbox{\rm if}\quad p\leq 2, \end{array}\right.
\end{equation}
\nd  where $x,y\in {\bf R}^N$ and  $(.,.)$ is the usual inner product of  ${\bf R}^N$.
\end{lemma}
\vskip.2cm

\nd The Hardy Inequality (cf. Br\'ezis \cite{B}) is:

\begin{theorem}\label{hardy}
There is a positive constant $C$ such that
$$
\int_{\Omega} \big|\frac{u}{d} \big|^{\beta} dx \leq  C \int_{\Omega} |\nabla u|^p,~ u\in \w.
$$
\end{theorem}

\noindent {\large \bf Proof of  lemma \ref{lem 3.3}}~ By the H\" older inequality,
\begin{equation}\label{eq B.4}
 \displaystyle\int_\Omega|\nabla u|^{p-1}|\nabla v|dx\\
                                \leq  ||u||_{1,p^{\prime}}||v||_{1,p},\end{equation}
\nd where $1/p + 1/p{^{\prime}} = 1$, and so the  expression
\begin{equation}\label{eq 3.9}
\langle-\Delta_p u,v\rangle :=\int_\Omega|\nabla u|^{p-2}\nabla u \cdot \nabla v dx, \quad u,v\in \w,
\end{equation}
\nd defines a  continuous, bounded (nonlinear)  operator  namely
$$\begin{array}{rlcrc}
\Delta_{p}  & : & \w &\longrightarrow  &W^{{-1},p'}(\Omega)\\
            &   & u  & \longmapsto     &\Delta_p u.  \end{array}
$$
\nd By (\ref{simon}), $-\Delta_p$  it is strictly monotone and coercive, that is
$$
\langle-\Delta_p u-(-\Delta_p v),u-v\rangle>0,~~  u,v\in\w,~ u\neq v
$$
\nd and
$$
\frac{\langle-\Delta_p u,u\rangle}{||u||_{1,p}}~ \stackrel{||u||_{1,p} \to \infty} {\longrightarrow}~\infty.
$$
\nd By the Browder-Minty Theorem,  $\Delta_{p}  :  \w \longrightarrow  W^{{-1},p'}(\Omega)$ is a homeomorphism. 
\vskip.1cm

\nd Consider 
 $$
 F_g(u)=\int_\Omega g u dx,~ u \in \w.
 $$
\nd {\bf Claim}~~ $F_g \in W^{-1,p'}(\Omega)$.
\vskip.2cm

\nd Assume for a while the {\bf Claim} has been proved. Since $-\Delta_p : \w \to W^{{-1},p'}(\Omega)$
is a homeomorphism, there is an only $u \in \w$ such that
$$
-\Delta_p u = F_g,
$$
\nd that is
$$
\langle-\Delta_p u,v\rangle = \int_\Omega gv dx,~ v\in \w
$$
\nd {\bf Verification of the Claim.} Let $V$ be an open neighborhood of $\partial \Omega$ such that
$0 < d(x)<1~\mbox{for}~ x\in V$ so that
$$
1<\frac{1}{d(x)^\beta}<\frac{1}{d(x)},~~ x\in V.
$$
\nd Now, if $v \in \w$ we have
$$
\displaystyle |F_g(v)|  \leq \displaystyle \int_\Omega|g||v|dx \\
= \displaystyle\int_{V^{c}}|g||v|dx +\int_{V}|g||v|dx\\
\leq \displaystyle C||v||_{1,p}+\int_{\Omega}\big|\frac{v}{d}\big|dx.
$$
\nd Applying the Hardy Inequality in the last term  above we get to,
$$
\displaystyle |F_g(v)| \leq  C ||v||_{1,p},
$$
\nd showing that  $F_g \in W^{-1,p'}(\Omega)$, proving the {\bf Claim}.
\vskip.5cm

\nd {\bf Regularity of $u$}: At first we treat the case $p=2$. By \cite{CRT} there is a solution $v$ of
$$
\left\{\begin{array}{rcll}\label{eq 1.1}
-\Delta v &=& \frac{1}{v^\beta} & \mbox{in}\  \Omega,\\
v &>&0 &\mbox{in}\  \Omega,\\
 v  &=& 0                                 &\mbox{on}\  \partial\Omega,
\end{array}\right.
$$
\nd which  belongs to $C^{1}(\overline{\Omega})$ and by the Hopf theorem $\frac{\partial v}{\partial \nu} < 0~\mbox{on}~\partial \Omega$. Since also $d \in C^{1}(\overline{\Omega})$ and $\frac{\partial d}{\partial \nu} < 0~\mbox{on}~\partial \Omega$  there a constant $C > 0$ such that
$$
v \leq C d~~\mbox{in}~~\Omega.
$$
\nd Moreover,
$$
-\Delta  v =\frac{1}{v^\beta}\geq \frac{C}{d^\beta}.
$$
\nd Consider the problem
$$\left\{\begin{array}{rrllr}
-\Delta \widetilde{u} & = & |g| & \mbox{in}& \Omega,\\
\widetilde{u}          & = & 0                    & \mbox{on} & \partial\Omega.
\end{array}\right.
$$
\nd By \cite[theorem B.1]{JST},
$$\widetilde{u}\in C^{1,\alpha}(\overline{\Omega})~ \mbox{and}~ ||\widetilde{u}||_{C^{1,\alpha}(\overline{\Omega})}\leq M_0.$$
\nd for some positive constant $M_0$. By the Maximum Principle,
$$
\widetilde{u} \leq  v \leq C d~\mbox{in}~\Omega.
$$
\nd Setting  $\overline{u} = u+\widetilde{u}$ we get
$$
-\Delta \overline{u} = g + |g|\geq0~ \mbox{in}~ \Omega
$$
\nd and by the arguments above,  $\overline{u}\leq C d~\mbox{in}~\Omega$. Thus,  as a consequence of \cite[theorem  B.1]{JST}, the are $\alpha \in (0,1)$ and $M_0 > 0$ such that
$$
\overline{u},\widetilde{u}\in C^{1,\alpha}(\overline{\Omega})~ \mbox{and}~ ||\overline{u}||_{C^{1,\alpha}(\overline{\Omega})},~ ||\widetilde{u}||_{C^{1,\alpha}(\overline{\Omega})}\leq M_0,
$$
\nd ending the proof of lemma  \ref{lem 3.3}    in the case $p = 2$.
\vskip.2cm

\nd In what follows we treat the  case $p > 1$. Let $u$ be a solution of (\ref{eq 3.8}). It follows that
$$
-\Delta_p u = g\leq\frac{C}{d^\beta}~~ \mbox{and}~ -\Delta_p (-u)=(-1)^{p-1}g\leq\frac{C}{d^\beta}.
$$
\nd By lemma \ref{lem 3.1}  the problem
$$
\left\{\begin{array}{rclr}
-\Delta_{p} v&=& \frac{C}{v^\beta}& \mbox{in}~ \Omega\\
v&=&0&  \mbox{on}~ \partial {\Omega}
\end{array}
\right.
$$
\nd admits an only positive solution  $v\in \w\cap C^{1,\alpha}(\overline\Omega)$ for some $\alpha\in(0,1)$ with $v \leq C d~\mbox{in}~\Omega$.
Hence,
$$
-\Delta_p(v)=\frac{C}{v^\beta}\geq \frac{1}{d^\beta}~ \mbox{in}~ \Omega.
$$
\nd Therefore,
$$
-\Delta_p|u|\leq \frac{C}{d^\beta}\leq -\Delta_p  v.
$$
\nd By the weak comparison principle,
$$
|u|\leq  v\leq C d~ \mbox{in}~ \Omega,
$$
\nd showing that  $u\in L^\infty(\Omega)$.  Pick $w\in C^{1,\alpha}(\overline{\Omega})$ such that
$$
-\Delta w=g~ \mbox{in}~ \Omega,~~  w=0~~  \partial\Omega.
$$
\nd We have
$$
div(|\nabla u|^{p-2}\nabla u-\nabla w)=0 \quad \mbox{em}\quad \Omega
$$
\nd in the weak sense. By  Lieberman \cite[theorem 1]{L} the proof of lemma \ref{lem 3.3}  ends.  $\hfill{\rule{2mm}{2mm}}$
\vskip.2cm

\noindent {\large \bf Proof of Corollary \ref{cor 3.2}}~ Existence of $u_\epsilon$ follows directly by lemma \ref{lem 3.3}. Moreover there are $M>0$ and $\alpha\in(0,1)$ such that
$$
||u||_{C^{1,\alpha}(\overline{\Omega})},~~||u_\epsilon||_{C^{1,\alpha}(\overline{\Omega})}<M.
$$
\nd By  V\'azquez \cite[theorem 5]{V},  $\frac{\partial u}{\partial\nu}<0~\mbox{on}~ \partial\Omega$ and recalling that $d \in C^{1}(\overline{\Omega})$ and $\frac{\partial d}{\partial \nu} < 0~\mbox{on}~\partial \Omega$   it follows that
\begin{equation}\label{eq 3.12}
 u \geq  C d~~ \mbox{in}~~ \Omega.
 \end{equation}
\nd Multiplying the equation
$$
 -\Delta_pu-(-\Delta_pu_\epsilon)  = g-\big( h\chi_{[d(x)>\epsilon]}+\widetilde{g}\chi_{[d(x)<\epsilon]}\big)
 $$
\nd by  $u-u_\epsilon$ and integrating we have
$$
\int_\Omega(|\nabla u|^{p-2}\nabla u-|\nabla u_\epsilon|^{p-2}\nabla u_\epsilon).\nabla (u-u_\epsilon)dx\leq2M\int_{d(x)<\epsilon}|g-\widetilde{g} |dx.
$$
\nd Using lemma  \ref{lem A.1}, we infer that $||u-u_\epsilon||_{1,p}\to 0$ as $\epsilon\to 0$. By the compact embedding  $C^{1,\alpha}(\overline{\Omega}) \hookrightarrow C^1(\overline{\Omega})$ it follows that
$$
||u-u_\epsilon||_{C^1(\overline{\Omega})}\leq \frac{C}{2}d,
$$
\nd and using (\ref{eq 3.12}),
$$
u_\epsilon \geq u-\frac{C}{2}d~ \geq~ u-\frac{u}{2} = \frac{u}{2}.~~~~~~~~~~~~~~~~  \hfill{\rule{2mm}{2mm}}
$$

\end{document}